\input amstex
\documentstyle{amsppt}
 \magnification=\magstep 1
\catcode`\@=11
\def\nologo{\let\logo@=\relax}
\catcode`\@=\active
\nologo

\vsize 7in
\topmatter

\title Socle degrees of Frobenius powers \endtitle 
  \leftheadtext{Kustin and Vraciu}
\rightheadtext{Socle degrees}
 \author Andrew R. Kustin and Adela N.  Vraciu\endauthor
 \address
Mathematics Department,
University of South Carolina,
Columbia, SC 29208\endaddress
\email kustin\@math.sc.edu vraciu\@math.sc.edu\endemail
\dedicatory In honor of Phillip  Griffith, on the occasion of his retirement.\enddedicatory
\abstract
Let $k$ be a field of positive characteristic $p$, $R$ be a Gorenstein graded $k$-algebra, and $S=R/J$ be an artinian quotient of $R$ by a homogeneous ideal.
 We ask how the socle degrees of $S$ are related to the socle degrees of $F_R^e(S)=R/J^{[q]}$. If $S$ has finite projective dimension as an $R$-module, then the socles of $S$ and $F_R^e(S)$ have the same dimension and the socle degrees are related by the formula:
$$D_i=qd_i-(q-1)a(R),$$
where
$$ d_1\le \dots\le d_{\ell}\quad\text{and}\quad D_1\le \dots \le D_{\ell}$$ are the socle degrees $S$ and $F_R^e(S)$, respectively, and   $a(R)$ is the $a$-invariant of the graded ring $R$, as introduced by Goto and Watanabe. 
We prove the converse when $R$ is a complete intersection. 
\endabstract
\endtopmatter
\document

Let $(R, \frak m)$ be a Noetherian graded algebra over a field of positive characteristic $p$, with irrelevant ideal $\frak m$.
We usually let $R=P/C$ with $P$ a polynomial ring, and $C$ a homogeneous ideal.
Let $J$ be an $\frak m$-primary homogeneous ideal in $R$. Recall that if $q=p^e$, then the $e^{\text{th}}$ Frobenius power of $J$ is the ideal $J^{[q]}$ generated by all $i^q$ with $i\in J$.
The basic question is:

\proclaim{Question}
How do the degrees of the minimal generators of $(J^{[q]} \:\!\frak m )/J^{[q]}$  vary with  $q$? 
\endproclaim

The largest of the degrees of a generators of the socle $(J\:\!\frak m)/J$ will be called the {\it top socle degree} of $R/J$.
The question of finding a linear bound for the top socle degree of $R/J^{[q]}$ has been considered by Brenner in \cite{2} from a different point of view; his main motivation there is finding inclusion-exclusion criteria for tight closure.

The answer to the   Question is well-known (although not explicitly stated in existing literature) in the case when $J$ has finite projective dimension; see Observation 1.7. 
We prove that the converse holds when $R=P/C$ is a complete intersection.
\proclaim{Theorem}Let $k$ be a field of positive characteristic $p$, $q=p^e$ for some positive integer $e$, $P$ be a positively graded polynomial ring over $k$, and   $R=P/C$ be a complete intersection ring with  $C$ generated by a homogeneous regular sequence.  Let $\frak m$ be the maximal homogeneous ideal of $R$, $J$ be a homogeneous $\frak m$-primary ideal in $R$, and $I$ be a lifting of $J$ to $P$. Let ${\ell}$ be the dimension of the socle $(J\:\!\frak m)/J$ of $R/J$ and  $d_1, \ldots, d_{\ell}$ be the degrees of the generators of the socle. Then the following statements are equivalent:

\medskip
\flushpar{\bf (a)} $\operatorname{pd}_R R/J< \infty$,

\medskip
\flushpar{\bf (b)} the socle $(J^{[q]}\:\!\frak m) /J^{[q]}$  of $R/J^{[q]}$ has dimension ${\ell}$ and the degrees of the generators are $qd_i-(q-1)a$, for $1\le i\le \ell$,  where $a$ denotes the $a$-invariant of $R$,

\medskip
\flushpar{\bf (c)} $(C+I)^{[q]}\:\! (C^{[q]}\:\! C)=C+I^{[q]}$, and

\medskip
\flushpar{\bf (d)} $I^{[q]}\cap C= (I\cap C)^{[q]}+CI^{[q]}$.

\endproclaim

\flushpar Of course, the general question 
\proclaim{Question}
How do the socle degrees of Frobenius powers $J^{[q]}$ encode homological information about the ideals $J^{[q]}$?
\endproclaim
\flushpar remains wide open and very compelling. 

\bigpagebreak

\flushpar{\bf 1.\quad  Preliminary notions.}

\medskip
 Let $k$ be a field of positive characteristic $p$. We say that the ring $R$ is a {\bf graded $k$-algebra} if $$\matrix\format\l\\\text{$R$ is non-negatively graded, $R_0=k$, and $R$ is finitely generated}\\\text{as a ring over $k$.}\endmatrix\tag 1.1$$
Every ring that we study in this paper is a graded $k$-algebra. In particular, ``Let $P$ be a polynomial ring'' means $P=k[x_1,\dots,x_n]$, for some $n$, and each variable has positive degree. Every calculation in this paper is homogeneous: all elements and ideals that we consider are homogeneous, all ring or module homomorphisms that we consider are homogeneous of degree zero. If $r$ is a homogeneous element of the ring $R$, then $|r|$ is the degree of $r$.  
The graded $k$-algebra $R$ has a unique homogeneous maximal ideal $$\frak m=\frak m_R=R_+=\bigoplus\limits_{i>0}R_i;$$ 
furthermore, $R$ has a unique graded canonical module $K_R$, which is equal to the graded dual of the  graded local cohomology module $\underline {\operatorname{H}}^n_{\frak m}(R)$,  where $n$ is the Krull dimension of $R$; that is,
$$K_R=\underline{\operatorname{Hom}}_R(\underline{\operatorname{H}}_\frak m^n(R),E_R),$$ 
for $E_R=\underline{\operatorname{Hom}}_k(R,R/\frak m)$  the injective envelope of $R/\frak m$ as a graded $R$-module. 
 (See, for example \cite{5, Def. 2.1.2}.) The $a$-invariant of $R$ is defined to be
$$a(R)=-\min\{m\mid (K_R)_m\neq 0\}=\max\{m\mid (\underline{\operatorname{H}}_{ \frak m}^n(R))_{m} \neq 0\}.$$
The definition of the $a$-invariant is rigged so that if $R$ is a Gorenstein graded $k$-algebra, then $K_R=R(a(R))$. 
When the ring $R$ is Cohen-Macaulay, there are many ways to compute $a(R)$. 
The main tool for these calculations, Proposition 1.2 below, may be found as Proposition 2.2.9 in \cite{5} or Proposition 3.6.12 in \cite{3}. The functor $\underline{\operatorname{Ext}}^c$ is underlined, as in \cite{5}, to emphasize the graded nature of this functor. All of our homological calculations are homogeneous; however, we will not take the trouble to underline every graded object. 
\proclaim{Proposition 1.2} If $R\to S$ is a graded surjection of graded $k$-algebras, and $R$ is Cohen-Macaulay, then 
$$K_S=\underline{\operatorname{Ext}}^c_R(S,K_R),$$ where $c=\dim R-\dim S$. In particular, if $S=R/C$ and the ideal $C$   is generated by the homogeneous regular sequence $f_1,\dots f_c$,
then $$\tsize K_{R/C}=\frac{K_R}{CK_R}(\sum\limits_{i=1}^c|f_i|).$$\endproclaim
\proclaim{Corollary 1.3}
\flushpar {\bf (a)} If $P$ is the polynomial ring $k[x_1,\dots,x_n]$,  then $a(P)=-\sum\limits_{i=1}^n|x_i|$.
\medskip\flushpar {\bf (b)} If $R$ is the complete intersection ring $P/C$ where $P$ is the polynomial ring $k[x_1,\dots,x_n]$ and $C$ is the ideal  in $P$  generated by the homogeneous regular sequence $f_1,\dots f_c$, then $a(R)= \sum\limits_{i=1}^c |f_i|-\sum\limits_{i=1}^n |x_i|$. 
\medskip\flushpar {\bf (c)} If $R\to S$ is a surjection of graded Cohen-Macaulay $k$-algebras, and  $S$ has finite projective dimension as an $R$-module, then 
$a(S)=a(R)+N$, where $N$ is the largest back twist in the minimal homogeneous resolution of $S$ by free $R$-modules. In other words, if 
$$0\to \bigoplus\limits_i R(-b_{c,i})\to \dots\to \bigoplus\limits_i R(-b_{1,i}) \to R \to S\to 0$$ is the minimal homogeneous resolution of $S$ by free $R$-modules, then  $N=\max\limits_{i}\{b_{c,i}\}$. 
\endproclaim

\definition{Definition} If $S$ is an artinian graded $k$-algebra, then the socle of $S$,
$$\operatorname{soc} S=0\:\! \frak m_S=\{s\in S\mid s\frak m_S=0\},$$ is a finite dimensional graded $k$-vector space:
$\operatorname{soc} S=\bigoplus\limits_{i=1}^{\ell}k(-d_i)$. We refer to the numbers $d_1\le d_2\le \dots\le d_{\ell}$ as the {\bf socle degrees of $S$}. \enddefinition

\proclaim{Observation 1.4} Let $R$ be an artinian Gorenstein graded $k$-algebra with socle degree $\delta$, and $J$ be a homogeneous ideal of $R$. If the socle degrees of $\frac RJ$ are $\{d_i\}$, then the minimal generators of $\operatorname{ann} J$ have degree $\{\delta-d_i\}$. \endproclaim
\demo{Proof}
Choose minimal generators $g_1, \ldots, g_s$ of $\operatorname{ann} J$. Gorenstein duality implies that $$\operatorname{ann}(g_1, \ldots, \hat{g_i}, \ldots, g_s) \not \subseteq \operatorname{ann}(g_i);$$ and thus, for each $i$, we can choose an  element $u_i \in \operatorname{ann}(g_1, \ldots, \hat{g_i}, \ldots, g_s)$, which represents a generator for the socle of $R/\operatorname{ann}(g_i)$. The ideals $J$ and  $\operatorname{ann}(g_1, \ldots, g_s)$ are equal and the socle of $R/\operatorname{ann}(g_1, \ldots, g_s)$ is minimally generated   by $u_1, \ldots, u_s$.
On the other hand, $u_i g_i$ generates the socle of $R$, so the degree of $u_i$ is equal to $\delta - |g_i|$. \qed
\enddemo

\proclaim{Proposition  1.5} If $S$ is an artinian graded $k$-algebra and $d_1\le \dots\le d_{\ell}$ are the socle degrees of $S$, then  the minimal generators of the canonical module $K_S$ have degrees $-d_{\ell}\le \dots\le -d_1$. \endproclaim
\demo{Proof} Let $P=k[x_1,\dots,x_n]$ be a polynomial ring which maps onto $S$. One may compute the degrees of the generators of $K_S$ as well the socle degrees of  $S$ in terms of the back twists in the minimal homogeneous resolution of $S$ as a $P$-module:
$$0\to \bigoplus_i P(-b_{n,i})\to \dots \to \bigoplus_i P(-b_{1,i})\to P \to S\to 0.$$ The canonical module $K_S$ is equal to $\underline{\operatorname{Ext}}_P^n(S,K_P)$, where $K_P=P(a(P))$ and $a(P)=-\sum\limits_{i=1}^n|x_i|$. It follows that the minimal homogeneous resolution of $K_S$ 
is $$0\to P(a(P))\to \dots \to  \bigoplus_i P(a(P)+b_{n,i})\to K_S\to 0;$$ and therefore, the minimal generators of $K_S$ (over either $S$ or $P$) have degrees $\{-a(P)-b_{n,i}\}$.  On the other hand, one may compute $\operatorname{Tor}_n^P(S,P/\frak m_P)$ in each coordinate (see for example \cite{6, Lemma 1.3}) in order to conclude that
$$\bigoplus_{i} k(-b_{n,i})= \operatorname{Tor}_n^P(S,k)=\operatorname{soc} S(a(P)).$$Thus, the socle degrees of $S$ are  equal to $\{a(P)+b_{n,i}\}$. \qed
\enddemo

\proclaim{Corollary 1.6} Let  $R\to S$ be a surjection of graded $k$-algebras with 
$S$ artinian, and $R$ Gorenstein. If  $\operatorname{pd}_RS$ finite, then the socle degrees of $S$ are
$\{b_i+a(R)\}$, where  the back twists in the minimal  homogeneous resolution of $S$ by free $R$-modules are $\{b_i\}$. \endproclaim
\demo{Proof} We know, from Proposition 1.2, that $K_S=\underline{\operatorname{Ext}}_R^{\dim R}(S,K_R)$, with $K_R=R(a(R))$; and therefore, 
$$0\to R(a(R))\to \dots\to \bigoplus\limits_iR(a(R)+b_i)\to K_S\to 0$$ is a minimal resolution of $K_S$ and the minimal generators of $K_S$ as an $R$-module, or as an $S$-module, have degrees $\{-a(R)-b_i\}$. Apply Proposition 1.5. \qed \enddemo

Let $R$ be a graded $k$-algebra. We write $^{\phi_R^e}\!R$ to represent the ring $R$ endowed with an $R$-module structure given by the $e^{\text{th}}$ iteration of the Frobenius endomorphism $\phi_R\: R\to R$. (If $r$ is a scalar in $R$ and $s$ is a ring element in  $^{\phi_R^e}\!R$, then $r\cdot s$ is equal to $r^qs\in {}^{\phi_R^e}\!R$, for $q=p^e$.) The Frobenius functor $F_R^e(\underline{\phantom{X}})=\underline{\phantom{X}}\otimes_R{}^{\phi_R^e}\!R$ is base change along the homomorphism $\phi_R^e$.  If $\pmb g$ is a matrix with entries in $R$, then $\pmb g^{[q]}$ is the matrix in which each entry of $\pmb g$ is raised to the power $q$. 
If $\Bbb G_1$ is the free module $\bigoplus_{i=1}^mR(-b_i)$, then $\Bbb G_1^{[q]}$ is the free module $\bigoplus_{i=1}^mR(-qb_i)$.   If $\pmb g$ is a matrix and $\pmb g\: \Bbb G_2\to \Bbb G_1$ is a map of free $R$-modules, then  $$\pmb g^{[q]}\: \Bbb G_2^{[q]}@> >> \Bbb G_1^{[q]}$$ is a very clean way to write  
$\left(\pmb g\: \Bbb G_2@> >> \Bbb G_1\right)\otimes_R{}^{\phi_R^e}\!R$. 
If $J$ is the of $R$-ideal $(a_1,\dots,a_m)$, then $J^{[q]}$ is the $R$-ideal $(a_1^q,\dots,a_m^q)$. The Frobenius functor is always right exact; so, in particular, 
$F_R^e(R/J)=R/J^{[q]}$. 

  \proclaim{Observation 1.7} Let $k$ be a field of positive characteristic $p$, $R\to S$ be  a  surjection of graded $k$-algebras  in the sense of $($1.1$)$, with $R$ Gorenstein  and $S$ artinian.
 If $S$ has finite projective dimension as an $R$-module, then the socles of $S$ and $F_R^e(S)$ have the same dimension; furthermore, if the socle degrees of 
$S$ and $F_R^e(S)$ are given by 
$$d_1\le d_2\le\cdots \le d_{\ell} \quad\text{and}\quad D_1\le D_2\le \cdots\le D_{\ell},$$ respectively, then
$$D_i=qd_i-(q-1)a(R),$$ for all $i$.\endproclaim
\demo{Proof}Consider the minimal homogeneous resolution $\Bbb F$ of $S$ by free $R$-modules. 
We know, from  the Theorem of Peskine and Szpiro \cite{7, Theorem 1.7}  that 
$F_R^e(\Bbb F)=\Bbb F^{[q]}$
is the minimal homogeneous  resolution of $F_R^e(S)$. If back twists of $\Bbb F$ are $\{b_i\mid 1\le i\le L\}$, then the back twists of $\Bbb F^{[q]}$ are $\{qb_i\}$. Use Corollary 1.6 to see that $L=\ell$, $d_i=b_i+a(R)$, and $D_i=qb_i+a(R)$, for all $i$. \qed
\enddemo

We prove the converse of Observation 1.7 under the assumption that $R$ is a complete intersection. Our main result is the following statement. 
\proclaim{Theorem 1.8}  Let $k$ be a field of positive characteristic $p$, $R\to S$ be  a  surjection of graded $k$-algebras  in the sense of $($1.1$)$, with $R$ a complete intersection  and $S$ artinian. Let $e$ be a positive integer, $q=p^e$, and $d_1\le \dots\le d_{\ell}$ be the socle degrees of $S$. If the socle of $F_R^e(S)$ has the same dimension as the socle of $S$, and the socle degrees of $F_R^e(S)$ are given by $D_1\le D_2\le \cdots\le D_{\ell}$, with $$D_i=qd_i-(q-1)a(R),$$ for all $i$,
then 
$$\operatorname{Tor}_1^R(S,{}^{\phi_R^e}\!R)=0.$$\endproclaim

\definition{The plan of attack 1.9}We express $R=P/C$, where $P$ is the polynomial ring $P=k[x_1,\dots,x_n]$, each variable has positive degree, and $C$ is a homogeneous Gorenstein ideal in $P$ of grade $c$. Let $I$ be a homogeneous $\frak m_P$-primary ideal in $P$ with $S=R/IR$. Let $T=P/I$. 

In Corollary 2.2, we convert  numerical information about the socle degrees of $S$ and $F_R^e(S)$ into numerical information about $\operatorname{Tor}_c^P(K_T,R)$ and $\operatorname{Tor}_c^P(K_{F_P^e(T)},R)$. In Proposition 3.1, the numerical information about $\operatorname{Tor}_c$'s is converted into the statement $$\operatorname{Tor}_1^R(M\otimes_PR,{}^{\phi_R^e}\!R)=0,$$ where $M$ is the the $(c-1)$-syzygy of the $P$-module $K_T$. This homological statement is expressed as a statement about ideals:
$$(C^{[q]}+I^{[q]})\:\!(C^{[q]}\:\! C)= C+I^{[q]}$$ in Proposition 4.1. 
In Proposition 5.1  we deduce 
$$I^{[q]}\cap C=(I\cap C)^{[q]}+CI^{[q]}.$$  This result is equivalent to 
$$\operatorname{Tor}_1^R(S,{}^{\phi_R^e}\!R)=0,$$as is recorded in Proposition 2.4. 
\enddefinition

We would like to prove that Theorem 1.8 continues to hold 
after one replaces the hypothesis that $R$ is a complete intersection with merely the hypothesis that $R$ is Gorenstein.  Three of our five steps (2.2, 4.1, and 2.4) work when $R$ is Gorenstein.
The arguments that we use in the other two steps (3.1 and 5.1) require that $R$ be a complete intersection; although in Proposition 6.2 we prove the ideal theoretic version of (3.1) under the hypothesis that $R$ is Gorenstein and F-pure. At any rate, if $R$ is a complete intersection and the conclusion of Theorem 1.8 holds, then the Theorem of Avramov and Miller \cite{1} (see also \cite{4}) guarantees that $S$ has finite projective dimension as an $R$-module. We are very curious to know if some form of  the Avramov-Miller result:
$$\operatorname{Tor}_1^R(M,{}^{\phi_R^e}\!R)=0\implies \operatorname{pd}_RM<\infty,$$for finitely generated $R$-modules $M$,  can be proven when $R$ is Gorenstein, but not necessarily a complete intersection. 

\bigpagebreak

\flushpar{\bf 2.\quad  Convert socle degrees into degrees of generators of $\operatorname{Tor}$-modules.}

\medskip
\proclaim{Lemma 2.1}Adopt the notation of 1.9. If the socle degrees of $S$ are $$\{d_i\mid 1\le i\le \ell\},$$ then the minimal generators of $\operatorname{Tor}_c^P(K_T(-a(P)),R)$ have degrees $$\{a(R)-d_i\mid 1\le i\le \ell\}.$$ 
\endproclaim

\demo{Proof}  Let $\Bbb G$ be the minimal homogeneous resolution of $R$ by free $P$-modules. Corollary 1.3~(c) tells us that $\Bbb G_c=P(a(P)-a(R))$. It follows that
$$\align &\operatorname{Tor}_c^P(K_T(-a(P)),R)=\operatorname{H}_c(K_T(-a(P))\otimes_P \Bbb G)\\&=\{\alpha\in K_T(-a(R))\mid C\alpha=0\}=\operatorname{Hom}_P(R,K_T(-a(R)).\endalign $$On the other hand, we have a surjection $T\to S$; so Proposition 1.2 guarantees 
$$K_S=\operatorname{Hom}_T(S,K_T)=\operatorname{Hom}_P(R,K_T).$$Thus, 
$$K_S(-a(R))=\operatorname{Hom}_P(R,K_T(-a(R)))=\operatorname{Tor}_c^P(K_T(-a(P)),R).$$ Apply Proposition 1.5. \qed\enddemo
 Lemma 2.1 also applies when the ideal $I$ is replaced by the ideal $I^{[q]}$; and consequently, 
if the socle degrees of $F_R^e(S)$ are $\{D_i\mid 1\le i\le L\}$, then the minimal generators of $\operatorname{Tor}_c^P(K_{F_P^e(T)}(-a(P)),R)$ have degrees $\{a(R)-D_i\mid 1\le i\le L\}$. Our conversion is complete.

\proclaim{Corollary 2.2} Retain the notation of 1.9. Assume that the socles of $S$ and $F_R^e(S)$ have the same dimension. Let 
$$d_1\le \dots\le d_{\ell} \quad\text{and}\quad D_1\le \dots\le D_{\ell}$$ be the socle degrees of $S$ and $F_R^e(S)$, respectively; and $$\gamma_1\le \dots\le \gamma_{\ell} \quad\text{and}\quad \Gamma_1\le \dots\le  \Gamma_{\ell},$$ be the minimal generator degrees of 
$$\operatorname{Tor}^P_c(K_T(-a(P)),R) \quad\text{and}\quad \operatorname{Tor}^P_c(K_{F_P^e(T)}(-a(P)),R),$$ respectively. Then 
$$D_i=qd_i-(q-1)a(R)\ \text{for all $i$}\iff \Gamma_i=q\gamma_i\ \text{for all $i$}.$$
\endproclaim

\flushpar In the notation of 1.9, let $ \Bbb F$ be the minimal homogeneous resolution of $K_T(-a(P))$. The proof of Proposition 1.5 shows that $F_P^e( \Bbb F)= \Bbb F^{[q]}$ is the minimal homogeneous resolution of $K_{F_P^e(T)}(-a(P))$ (rather than some other shift of $K_{F_P^e(T)})$; hence,
$$\align F_P^e(K_T(-a(P)))&{}= K_{F_P^e(T)}(-a(P)),\\
 \operatorname{Tor}_i^P(K_T(-a(P)),R)&{}=\operatorname{H}_i( \Bbb F\otimes_PR),\ \text{and}\\
\operatorname{Tor}_i^P(F_P^e(K_T(-a(P))),R)&{}=\operatorname{H}_i( \Bbb F^{[q]}\otimes_PR),\endalign$$ for all $i$.
We focus on $i=c$. Let $M$ be the $(c-1)$-syzygy of $K_T(-a(P))$. The beginning of the minimal homogeneous resolution of $M$ is 
$$\dots\to  \Bbb F_{c+1}\to  \Bbb F_c\to  \Bbb F_{c-1}\to M\to 0.$$ One may calculate $\operatorname{Tor}_1^R(M\otimes_PR,{}^{\phi_R^e}\!R)$ as 
$$\frac{\operatorname{H}_c( \Bbb F^{[q]}\otimes_PR)}{\left(\operatorname{H}_c( \Bbb F\otimes_PR)\right)^{[q]}},$$as described in Observation 2.3. 
The conclusion of Corollary 2.2 tells us that if $[\bar z_1]{,\allowmathbreak }\dots{,\allowmathbreak } [\bar z_{\ell}]$ is a minimal generating set for 
$\operatorname{H}_c( \Bbb F\otimes_PR)$, then the elements $[\bar z_1^{[q]}],\dots,[\bar z_{\ell}^{[q]}]$, in $\operatorname{H}_c( \Bbb F^{[q]}\otimes_PR)$, have the correct degrees to be a minimal generating set for $\operatorname{H}_c( \Bbb F^{[q]}\otimes_PR)$. In Proposition 3.1 we prove that, when $R$ is a complete intersection, then  $[\bar z_1^{[q]}],\dots,[\bar z_{\ell}^{[q]}]$ do indeed generate $\operatorname{H}_c( \Bbb F^{[q]}\otimes_PR)$.

 \proclaim{Observation 2.3}Let $P\to R$ be a surjection of graded $k$-algebras, with $P$ a polynomial ring, and let $M$ be a finitely generated graded $P$-module. Then there is an exact sequence of graded $R$-modules:
$$F_R^e\left(\operatorname{Tor}_1^P(M,R)\right)\to \operatorname{Tor}_1^P(F_P^e(M),R) \to \operatorname{Tor}_1^R(M\otimes_PR,{}^{\phi_R^e}R)\to 0.$$
\endproclaim
\demo{Proof}Let $( \Bbb N,\pmb n)$ be the minimal homogeneous resolution of $M$ by free $P$-modules. The functor $F_P^e(\underline{\phantom{X}})$ is exact; so, $ \Bbb N^{[q]}$ is the minimal homogeneous resolution of $F_P^e(M)$ by free $P$-modules; and $\operatorname{Tor}_1^P(F_P^e(M),R)$ is equal to 
$$\operatorname{H}_1(F_P^e( \Bbb N)\otimes_PR)=\operatorname{H}_1(F_R^e( \Bbb N\otimes_PR)).$$ The functors $F_P^e(\underline{\phantom{X}})\otimes_PR$ and $F_R^e(\underline{\phantom{X}}\otimes_PR)$ are equal 
because the homomorphisms $$\CD P@>\phi_{P}^e >> P\\@V \text{quot. map} VV @V V\text{quot. map} V\\R@>\phi_{R}^e >> R\endCD$$ commute.
Let $\bar{\phantom{x}}$ denote the functor $\underline{\phantom{X}}\otimes_PR$.
Select elements
$z_1,\dots,z_{\ell}$  of $ \Bbb N_1$ so that  $\bar z_1,\dots,\bar z_{\ell}$ are cycles in $ \Bbb N\otimes_PR$ and the homology classes 
$[\bar z_1],\dots,[\bar z_{\ell}]$ are a minimal generating set for $\operatorname{H}_1( \Bbb N\otimes_PR)=\operatorname{Tor}_1^P(M,R)$. It is clear that $z_i^{[q]}$ is an element of $F_P^e( \Bbb N_1)$ with $\overline{z_i^{[q]}}=\bar z_i^{[q]}$ a cycle in $F_P^e( \Bbb N)\otimes_PR=F_R^e( \Bbb N\otimes_PR)$, for each $i$.

The technique of killing cycles tells us that  
$$\Bbb M\:\quad \bar  \Bbb N_2\oplus\bigoplus\limits_{i=1}^{\ell}R(-|z_i|)@> \bmatrix \bar {\pmb  n}_2&\bar z_1&\dots \bar z_{\ell}\endbmatrix>>\bar  \Bbb N_1@> \bar {\pmb  n}_1 >> \bar  \Bbb N_0\to \bar M\to 0$$ is the beginning of a homogeneous resolution of $\bar M$ by free $R$-modules. It follows that 
$$\operatorname{Tor}_1^R(\bar M,{}^{\phi_R^e}\!R)=\operatorname{H}_1(\Bbb M\otimes_R{}^{\phi_R^e}\!R)=\frac{\operatorname{H}_1(F_R^e(\bar  \Bbb N))}{([\bar z_1^{[q]}],\dots,[\bar z_{\ell}^{[q]}])}=\frac{\operatorname{Tor}_1^P(F_P^e(M),R)}{([\bar z_1^{[q]}],\dots,[\bar z_{\ell}^{[q]}])}. \qed $$
\enddemo

The following result is an application of the technique of Observation 2.3. 
\proclaim{Proposition 2.4}Let $R=P/C$ and $T=P/I$, where   $I$ and $C$ are homogeneous ideals in the polynomial ring $P$. Then $$\operatorname{Tor}_1^R(T\otimes_PR, {}^{\phi_R^e}R)=  \frac{I^{[q]}\cap C} {(I\cap C)^{[q]}+ I^{[q]}C}.$$
\endproclaim
\demo{Proof}Apply Observation 2.3 to see that $$\operatorname{Tor}_1^R(T\otimes_PR, {}^{\phi_R^e}R)=  \frac{\operatorname{H}_1(  \Bbb N^{[q]}\otimes_PR)}{([\bar z_1^{[q]}],\dots,[\bar z_{\ell}^{[q]}])},$$ where
$( \Bbb N,\pmb n)$ is the minimal homogeneous resolution of $T$ by free $P$-modules,   $\bar{\phantom{x}}$ is the functor $\underline{\phantom{X}}\otimes_PR$, and 
$z_1,\dots,z_{\ell}$  are elements of  of $ \Bbb N_1$ with  
$[\bar z_1],\dots,[\bar z_{\ell}]$  a minimal generating set for $$\tsize \operatorname{H}_1( \Bbb N\otimes_PR)=\operatorname{Tor}_1^{P}(\frac PI,\frac PC)=\frac{I\cap C}{IC}.$$ Observe that  $I\cap C=(\pmb n_1(z_1),\dots,\pmb n_1(z_{\ell}))+IC$.  Observe also, that
$$\tsize \operatorname{H}_1(  \Bbb N^{[q]}\otimes_PR)=\operatorname{Tor}_1^{P}(\frac P{I^{[q]}},\frac PC)=\frac{I^{[q]}\cap C}{I^{[q]}C}.$$The isomorphism $\operatorname{H}_1(  \Bbb N^{[q]}\otimes_PR)\to \frac{I^{[q]}\cap C}{I^{[q]}C}$ carries $[\bar z_i^{[q]}]$ to the class of $(\pmb n_1(z_i))^q$. \qed
\enddemo

\bigpagebreak

\flushpar{\bf 3.\quad  Degree considerations concerning $\operatorname{Tor}_1$.}
\medskip

Proposition 3.1 is a general statement that says that if the degrees of the minimal generators of 
$$\operatorname{Tor}_1^P(M,R)\quad\text{and}\quad \operatorname{Tor}_1^P(F_P^e(M),R)$$ are related in the appropriate manner, then $M\otimes_PR$ has finite projective dimension as an $R$-module. When the notation of 1.9 and the hypothesis of Theorem 1.8 are in effect, then, as we saw after Corollary 2.2,  Proposition 3.1  may be applied with $M$ equal to the $(c-1)$-syzygy of $K_T(-a(P))$. 
\proclaim{Proposition 3.1}
Let $P\to R$ be a surjection of graded $k$-algebras, with $P$ a polynomial ring and $R$ a complete intersection, and let $M$ be a finitely generated graded $P$-module.
Suppose that the minimal generators of $\operatorname{Tor}_1^P(M,R)$ have degrees $\{\gamma_i\mid 1\le i\le \ell\}$. If  the minimal generators of $\operatorname{Tor}_1^P(F_P^e(M),R)$ have  degrees $\{q\gamma_i\mid 1\le i\le \ell\}$, then $\operatorname{Tor}_1^R(M\otimes_PR,{}^{\phi_R^e}\!R)=0$. \endproclaim

\demo{Proof}Inflation of the base field $k\to K$ gives rise to faithfully flat extensions $P\to P\otimes_kK$ and $R\to R\otimes_kK$. Consequently,  we may assume  that $k$ is a perfect field. Let $C$ be the ideal in $P$ with $R=P/C$, and let
$f_1,\dots,f_c$ be a homogeneous regular sequence in $P$  that generates $C$. 
We retain the notation from the proof of Observation 2.3. So,
 $$\Bbb N\:\quad \Bbb N_2@> {\pmb n}_2>> \Bbb N_1 @> {\pmb n}_1>> \Bbb N_0\to M\to 0$$ is the beginning of the  minimal homogeneous  resolution of $M$ by free $P$-modules, $\bar{\phantom{x}}$ is the functor $\underline{\phantom{X}}\otimes_PR$,  and $z_1,\dots,z_{\ell}$ are elements of $\Bbb N_1$
with $[\bar z_1],\dots,[\bar z_{\ell}]$ a minimal generating set for $\operatorname{H}_1(\bar \Bbb N)=\operatorname{Tor}_1^P(M,R)$.
We will prove that $[\bar z_1^{[q]}], \dots, [\bar z_{\ell}^{[q]}]$ generates 
$$\operatorname{H}_1(\bar\Bbb N^{[q]})=\operatorname{Tor}_1^P(F_P^e(M),R).$$
 For each integer $\delta$, let $V_{\delta}$ be  the vector space  $$V_{\delta}=\frac{\operatorname{H}_1(\bar\Bbb N^{[q]})}{W_{\delta}'+W_{\delta}''},$$where  
 $W_{\delta}'$ is the $P$-submodule $$\sum\limits_{\delta<i}\left[\operatorname{H}_1(\bar\Bbb N^{[q]})\right]_i$$ of $\operatorname{H}_1(\bar\Bbb N^{[q]})$ and  $W_{\delta}''$ is the $P$-submodule of $\operatorname{H}_1(\bar\Bbb N^{[p]})$ generated by $$\sum\limits_{i<\delta}\left[\operatorname{H}_1(\bar\Bbb N^{[q]})\right]_i.$$
 Let $X$ be a homogeneous minimal generating set for the $P$-module  $\operatorname{H}_1(\bar\Bbb N^{[q]})$. Let $X_{\delta}$, be the subset of $X$  which consists of  those generators which have degree equal to $\delta$. Notice that the set  $X_{\delta}$ is a basis the vector space  $V_{\delta}$ over $k$.
Let $$Y_{\delta}=\{[\bar z_i^{[q]}]\in \operatorname{H}_1(\bar\Bbb N^{[q]})\mid 1\le i\le \ell\quad\text{and}\quad |z_i^{[q]}|=\delta\}.$$
Our hypothesis guarantees that the dimension of the vector space $V_{\delta}$  is exactly equal to the number of elements in $Y_{\delta}$. Once we prove that the elements of $Y_{\delta}$ are linearly independent in $V_{\delta}$, then we will know that $Y_{\delta}$ is a basis for $V_{\delta}$ and the elements of $Y_{\delta}$ are part of a minimal generating set of   $\operatorname{H}_1(\bar\Bbb N^{[q]})$. 

We work by induction on $\delta$ starting with $\delta=0$ and then looking at successively higher values of $\delta$.
When $\delta$ is small, then  $Y_{\delta}$ is the empty set,  $V_{\delta}$ is zero, and everything is fine. Now we work on the inductive step. There is nothing to do unless $\delta$ is a multiple of $q$. Relabel the $z_i$, if necessary,
and select the integer $\lambda$ so that $z_1,\dots,z_{\lambda}$ have degree less than $\delta/q$, and $z_{\lambda+1},\dots,z_{\ell}$ have degree at least  $\delta/q$. 
Consider a non-trivial $k$-linear combination of of the elements of $Y_{\delta}$. The field $k$ is closed under the taking of $q^{\text{th}}$ roots; so this linear combination is equal to $[\bar z^{[q]}]$ where $z$ is a non-trivial $k$-linear combination of those $z_i$ that have degree equal to $\delta/q$. So
$$\text{$[\bar z]$ is not zero in } \frac{\operatorname{H}_1(\bar\Bbb N)}{
P([\bar z_1],\dots, [\bar z_{\lambda}])}.\tag 3.2$$
We assume that  
$$[\bar z^{[q]}]=0\text{ in }\frac{\operatorname{H}_1(\bar\Bbb N^{[q]})}{
W_{\delta}''};\tag 3.3$$
and we will show that this assumption leads to a contradiction. Keep in mind that the induction hypothesis guarantees that
$$ W_{\delta}''
=
P([\bar z_1^{[q]}],\dots, [\bar z_{\lambda}^{[q]}]).\tag 3.4$$
We know that $\bar z$ and $\bar z_1,\dots, \bar z_{\lambda}$ all are cycles in $\bar \Bbb N$; so $${\pmb n}_1z\in C\Bbb N_0\quad\text{and}\quad {\pmb n}_1z_i\in C\Bbb N_0\text{ for $1\le i\le \lambda$.}\tag 3.5$$
We introduce a notational convenience. Let
$$\Bbb Q_2=\Bbb N_2\oplus\bigoplus\limits_{i=1}^{\lambda} P(-|z_i|)$$ and let $\pmb q_2\:\Bbb Q_2\to \Bbb N_1$ be the map
$$\pmb q_2=\bmatrix  \pmb n_2& z_1&\dots&z_{\lambda}\endbmatrix.$$ Of course, the map
$$\pmb q_2^{[q]}\:\Bbb Q_2^{[q]}\to \Bbb N^{[q]}_1$$ also now has meaning. Assumption 3.3, together with the induction hypothesis (3.4), 
tells us that 
$$z^{[q]}\in \operatorname{im} \pmb q_2^{[q]}+C\Bbb N^{[q]}_1;$$which is the base case for the following induction.
We prove that if ${1\le t\le c(q-1)}$, then
$$z^{[q]}\in \operatorname{im} \pmb q_2^{[q]}+C^{[q]}\Bbb N^{[q]}_1 +C^t\Bbb N^{[q]}_1\implies z^{[q]}\in \operatorname{im} \pmb q_2^{[q]}+C^{[q]}\Bbb N^{[q]}_1 +C^{t+1}\Bbb N^{[q]}_1.\tag 3.6$$
As soon as (3.6) is established, then 
$$z^{[q]}\in \operatorname{im}\pmb q_2^{[q]}+C^{[q]}\Bbb N^{[q]}_1$$
because $C^{c(q-1)+1}\subseteq C^{[q]}$. Now the proof is complete because we apply Observation 3.7 to 
$$z^{[q]}\in \operatorname{im} \bmatrix \pmb q_2 &f_1I&\cdots &f_cI\endbmatrix^{[q]}$$
to conclude that $$z\in \operatorname{im}\bmatrix \pmb q_2 &f_1I&\dots&f_cI\endbmatrix,$$ and this violates (3.2).

We prove (3.6).
For each $c$-tuple $\alpha=(\alpha_1,\dots,\alpha_c)$ of non-negative integers, with $\alpha_i<q$ for all $i$, and $\sum \alpha_i=t$, there exists $y_{\alpha}\in \Bbb N^{[q]}_1$ such that 
$$z^{[q]}-\sum\limits_{\alpha}f_1^{\alpha_1}\cdots f_c^{\alpha_c}y_{\alpha}\in C^{[q]}\Bbb N_1^{[q]}+\operatorname{im} \pmb q_2^{[q]}.$$
Fix a $c$-tuple $\alpha$. Apply ${\pmb n}_1^{[q]}$ and use (3.5) to see that $$f_1^{\alpha_1}\cdots f_c^{\alpha_c}{\pmb n}_1^{[q]}y_{\alpha}\in (f_1^{\alpha_1+1},\dots,f_c^{\alpha_c+1})\Bbb N^{[q]}_1.$$
It follows, from the fact that $f_1,\dots,f_c$ is a regular sequence, that $${\pmb n}_1^{[q]}y_{\alpha}\in C\Bbb N^{[q]}_1.$$
So, $\bar y_{\alpha}$ is a one-cycle, of degree less than $\delta$, in $\Bbb N^{[q]}$.
The induction hypothesis (3.4) tells us that
$$y_{\alpha}\in \operatorname{im}\pmb q_2^{[q]}+C\Bbb N_1^{[q]},$$ and (3.6) is established.
 \qed
\enddemo

We close this section with a quick application of the flatness of the Frobenius functor for regular rings. 

\proclaim{Observation 3.7}Let $P$ be a polynomial ring, $\pmb q\:\Bbb Q_2\to \Bbb Q_1$ be a homomorphism of graded $P$-modules, and $z$ be an element of $\Bbb Q_1$. If $z^{[q]}$ is in the image of $\pmb q^{[q]}$, then $z$ is in the image of $\pmb q$. \endproclaim

\demo{Proof}Let $\tilde \Bbb Q_2$ be the graded free module $\Bbb Q_2\oplus P(-|z|)$ and $\tilde {\pmb q}$ be the map of graded free modules
$$\tilde {\pmb q}=\bmatrix \pmb q& z\endbmatrix\: \tilde \Bbb Q_2\to \Bbb Q_1.$$
The hypothesis ensures the existence of $h\in \Bbb Q_2^{[q]}$ with  $\left[\smallmatrix -h\\1\endsmallmatrix\right]$ is in the kernel of $\tilde {\pmb q}^{[q]}$. The Frobenius functor $\underline{\phantom{X}}\otimes_P{}^{\phi_{P}^e}\!P$ is flat; so, there exist $\left[\smallmatrix t_i\\b_i\endsmallmatrix\right]$ in $\ker \tilde {\pmb q}$ and $a_i\in P$ such that $\sum a_i\left[\smallmatrix t_i\\b_i\endsmallmatrix\right]^{[q]}=\left[\smallmatrix -h\\1\endsmallmatrix\right]$. Degree considerations tell us that $b_i$ is a unit, for some $i$. For this $i$, we have 
${\pmb q} (\frac{-t_i}{b_i})=z$.
\qed \enddemo

\bigpagebreak

\flushpar{\bf 4.\quad  We interpret $\operatorname{Tor}_1$ of the syzygy in terms of ideals.}

\medskip
Recall, from the beginning of section 3, that if the notation of 1.9 and the hypotheses of Theorem 1.8 are in effect, then $\operatorname{Tor}_1^R(M\otimes_PR,{}^{\phi_R^e}\!R)=0$, where $M$ is the $(c-1)$-syzygy of $K_T(-a(P))$. In this section, we interpret this $\operatorname{Tor}$-module in terms of ideals. Our interpretation continues to hold even when the hypotheses of Theorem 1.8 are not in effect.

\proclaim{Proposition 4.1}
Adopt the notation of 1.9 and  let 
   $M$ be the $(c-1)$-syzygy in the minimal homogeneous resolution of $K_{T}(-a(P))$ by free $P$-modules. Then 
 $$\operatorname{Tor}_1^R(M\otimes_P R,{}^{\phi_R^e}R)=\operatorname{Hom}_{P}\left(\frac{(C^{[q]}+I^{[q]})\:\! (C^{[q]}\:\! C)}{C+I^{[q]}}, \frac{P}{A^{[q]}}(N)\right),\tag 4.2$$
where $A\subseteq I$ is a homogeneous $\frak m_P$-primary  Gorenstein ideal and $N$ is equal to $a(P/A^{[q]})-a(R)$; furthermore,
$$\operatorname{Tor}_1^R(M\otimes_P R,{}^{\phi_R^e}R)=0 \iff \frac{(C^{[q]}+I^{[q]})\:\! (C^{[q]}\:\! C)}{C+I^{[q]}}=0.$$
\endproclaim

\demo{Proof} In our argument,  we assume that $A\subseteq I$ have the same grade without assuming that these ideals are $\frak m_P$-primary, and we prove that
 $$\operatorname{Tor}_1^R(M\otimes_P R,{}^{\phi_R^e}R)=\frac {A^{[q]}\:\! (C+I^{[q]})}{(A\:\! (C+I))^{[q]}(C^{[q]}\:\! C)+ A^{[q]}}(N).\tag 4.3$$
Once (4.3) is established, then, when  $A$ is a $\frak m_P$-primary ideal,  Gorenstein duality may be employed to show that  to the right side of (4.3) and the right side of (4.2) are both equal 
to$$\frac {A^{[q]}\:\! (C+I^{[q]})}{A^{[q]}\:\! \left((C+I)^{[q]}\:\! (C^{[q]}\:\! C)\right)}(N)
.$$

  Let $(\Bbb G,\pmb g_{\bullet})$ be the minimal homogeneous resolution of $R$ by free $P$-modules and $\pmb y_{\bullet}\: \Bbb G^{[q]}\to \Bbb G$ be a map of complexes which lifts the natural quotient map $P/C^{[q]}\to R$. 
The ideal $C$ is Gorenstein of grade $c$; so   $ \Bbb G_c= P(a(P)-a(R))$.  
The map $\pmb y_c\:\Bbb G_c^{[q]}\to \Bbb G_c$ is multiplication by $y$, for some element $y$ in $P$ of degree $(q-1)(a(R)-a(P))$.  
 We know, from linkage theory,  that
$$C^{[q]}\:\! C=(y,C^{[q]}).$$
 Use the surjections $\frac PA\to T$ and $\frac P{A^{[q]}}\to F^e_P(T)$ to calculate
$$\tsize K_T=\frac {A\:\! I}A(a(P/A))\quad\text{and}\quad 
K_{F^e_P(T)}=\frac {A^{[q]}\:\! I^{[q]}}{A^{[q]}}(a(P/A^{[q]})).$$
It follows that   $\operatorname{H}_c(K_T(-a(P))\otimes_P \Bbb G)$ is equal to 
$$\tsize \{ \alpha\in \frac {A\:\! I}A(a(P/A)-a(R)) \mid \alpha C=0\}= \frac { A\:\! (I+ C)}A(a(P/A)-a(R)) $$and  $$\operatorname{H}_c(K_{F^e_P(T)}(-a(P))\otimes_P \Bbb G)= \frac { A^{[q]}\:\! (I^{[q]}+ C)}{A^{[q]}}(N) . $$

Let $(\Bbb F,\pmb f)$
be the minimal homogeneous  resolution of $K_T(-a(P))$ by free $P$ modules.
We saw in Observation 2.3 that 
$$\operatorname{Tor}_1^R(M\otimes_PR,{}^{\phi^e_R}\!R)=\frac{\operatorname{H}_c(\Bbb F^{[q]}\otimes_PR)}{([\bar z_1^{[q]}],\dots,[\bar z_{\ell}^{[q]}])},$$ where $z_1,\dots,z_{\ell}$ are elements in $\Bbb F_c$ with $[\bar z_1],\dots,[\bar z_{\ell}]$ a minimal generating set for $\operatorname{H}_c(\Bbb F\otimes_PR)$. 
We use the isomorphisms
$$
 \matrix \format\l&\l\\ \operatorname{H}_c(\Bbb F\otimes_P R)= \operatorname{H}_c(\operatorname{Tot}(\Bbb F\otimes_P\Bbb G))&{}= \operatorname{H}_c(K_T(-a(P))\otimes_P \Bbb G)\\\vspace{4pt} &{}= \frac{ A\:\! (C+I)}{A}(a(P/A)-a(R))\endmatrix  \tag 4.4$$
and 
$$\matrix \format\l&\l\\  \operatorname{H}_c(\Bbb F^{[q]}\otimes_P R)= \operatorname{H}_c(\operatorname{Tot}(\Bbb F^{[q]}\otimes_P\Bbb G))&{}= \operatorname{H}_c(K_{F_P(T)}(-a(P))\otimes_P \Bbb G)\\\vspace{4pt}&{}= 
\frac{A^{[q]}\:\! (C+I^{[q]})}{A^{[q]}}(N)\endmatrix\tag 4.5$$
to re-express $\operatorname{Tor}_1^R(M\otimes_P R,{}^{\phi^e_R}R)$ as a subquotient of  $P(N)$. 
First, we give an explicit description of  the isomorphism (4.4).
 Let $w_c$ be an element of $\Bbb F_c$ with $\pmb f_c(w_c)\in C\Bbb F_{c-1}$. Each row and column of the double complex $\Bbb F\otimes \Bbb G$ is exact; and therefore, for each $i$, with $0\le i\le c$, there exists $w_i\in \Bbb F_i\otimes \Bbb G_{c-i}$ such that $$(\pmb f_{i}\otimes 1)(w_i)=(1\otimes \pmb g_{c-i+1})(w_{i-1}).\tag 4.6$$ In particular, $w_0$ is an element of $\Bbb F_0\otimes \Bbb G_c$.  The isomorphism (4.4) sends the homology class  $[\bar w_c]$ to the image of   $w_0$ in $\frac{ A\:\! (C+I)}{A}(a(P/A)-a(R))$. 
In a similar manner, if $W_c$ is  an element of $\Bbb F_c^{[q]}$ with $\pmb f_c^{[q]}(W_c)\in C\Bbb F_{c-1}^{[q]}$, then the isomorphism (4.5)
sends the homology class $[\bar W_c]$ in $\operatorname{H}_c(\Bbb F^{[q]}\otimes_P R)$ to the image of  $W_0$ in $\frac{A^{[q]}\:\! (C+I^{[q]})}{A^{[q]}}(N)$, where $W_i\in \Bbb F_i^{[q]}\otimes \Bbb G_{c-i}$ and 
$$(\pmb f_{i}^{[q]}\otimes 1)(W_i)=(1\otimes \pmb g_{c-i+1})(W_{i-1}).\tag 4.7$$

We finish the argument by showing that the submodule $([\bar z_1^{[q]}],\dots,[\bar z_{\ell}^{[q]}])$ of $\operatorname{H}_c(\Bbb F^{[q]}\otimes_PR)$ is sent to the submodule 
$\frac{(A\:\! (C+I))^{[q]}y+A^{[q]}}{A^{[q]}}(N)$ of $\frac{A^{[q]}\:\! (C+I^{[q]})}{A^{[q]}}(N)$ under the isomorphism (4.5).
Let $\bar w_c$ be a cycle in $\Bbb F\otimes R$, for some element $w_c$   of $\Bbb F_c$.  We are given the family $\{w_i\}$ with $w_i\in \Bbb F_i\otimes \Bbb G_{c-i}$ such that (4.6) holds.  If $w_i=\sum\limits_j u_{i,j}\otimes v_{c-i,j}$ with $u_{i,j}\in \Bbb F_i$ and $v_{c-i,j}\in \Bbb G_{c-i}$, then let $$W_i= \sum\limits_j u_{i,j}^{[q]}\otimes \pmb y_{c-i}(v_{c-i,j}^{[q]})\in \Bbb F_i^{[q]}\otimes \Bbb G_{c-i}.$$ A short calculation shows that that (4.7) holds for $\{W_i\}$ and we conclude that if  $a$ in $\frac{ A\:\! (C+I)}{A}(a(P/A)-a(R))$ is  the image of the homology class $[\bar w_c]$ under the isomorphism (4.4), then $ya^q$ in $\frac{A^{[q]}\:\! (C+I^{[q]})}{A^{[q]}}(N)$ is the image of the homology class $[\bar w_c^{[q]}]$ under the isomorphism (4.5). \qed\enddemo

\bigpagebreak

\flushpar{\bf 5.\quad  The key calculation.}

\medskip
In section 3, we proved that if the socle hypothesis of Theorem 1.8 is in effect and $M$ is the $(c-1)$-syzygy of $K_T$ as a $P$-module, then $\operatorname{Tor}_1^R(M\otimes_PR,{}^{\phi_R^e}\!R)=0$. Our goal is to prove that
$\operatorname{Tor}_1^R(T\otimes_PR,{}^{\phi_R^e}\!R)=0$. Homological arguments in sections 2 and 4 connect these $\operatorname{Tor}$-modules to quotients of ideals. In the present section we show how information about the $\operatorname{Tor}$-module of $M$ gives information about the $\operatorname{Tor}$-module of $T$, when $R$ is a complete intersection.

\proclaim{Proposition 5.1}Let $P$ be a regular ring of positive characteristic $p$, and let  $C$ and $I$ be ideals in $P$.
Assume that  $C$ is generated by the regular sequence $f_1,\dots,f_c$ and that
$$(C+I)^{[q]}\:\! y= C+I^{[q]},$$where $y=(f_1\cdots f_{c})^{q-1}$. Then
$$I^{[q]}\cap C=(I\cap C)^{[q]}+CI^{[q]}.$$
\endproclaim
\demo{Proof}Notice that $C^{[q]}\:\! C=(y)+C^{[q]}$.  Take $\xi\in I^{[q]}\cap C$. We prove that
if $1\le t\le c(q-1)$, then
$$\xi\in C^t+C^{[q]}+CI^{[q]}\implies \xi\in C^{t+1}+C^{[q]}+CI^{[q]}.\tag 5.2$$
Of course, we know that the hypothesis of (5.2) holds for $t=1$.  Once we have established (5.2), then, since  $C^{c(q-1)+1}\subseteq C^{[q]}$, we know that
$$\xi\in I^{[q]}\cap C^{[q]}+CI^{[q]}=(I\cap C)^{[q]}+CI^{[q]},$$
because the Frobenius functor on $P$ is flat. Now we prove (5.2).
Write $\xi$ as an element of $C^{[q]}+CI^{[q]}$ plus
$$\sum_{\alpha}b_{\alpha}f_1^{\alpha_1}\cdots f_{c}^{\alpha_c},$$
where $\alpha=(\alpha_1,\dots,\alpha_c)$ varies over all $c$-tuples of non-negative integers with $\alpha_i<q$ for all $i$ and 
$\sum_{i=1}^c \alpha_i=t$. Fix an index $\alpha$. Observe that 
$$f_1^{q-1-\alpha_1}\cdots f_c^{q-1-\alpha_c}\xi$$
is equal to $b_{\alpha}y$ plus an element of $C^{[q]}+I^{[q]}$. The hypothesis tells us that $b_{\alpha}$ is in $C+I^{[q]}$; (5.2) is established, and the proof is complete. \qed\enddemo

\bigpagebreak

\flushpar{\bf 6.\quad   The Gorenstein F-pure case.}

\medskip
The question of whether the conclusion of Theorem 1.8 holds when $R$ is Gorenstein is still open. In this section, we include partial results in this direction. Recall that the ring $R$ of positive prime characteristic $p$  is F-pure 
  if whenever $J$ is an ideal of $R$ and $x$ is an element of $R$ with  $x \notin J$, then $x^q \notin J^{[q]}$ for all $q=p^e$.

First note that the top socle degree ($\operatorname{tsd}$) of a Frobenius power is always at least equal to the ``expected'' top socle degree:
\proclaim{Proposition 6.1}Let $k$ be a field of positive characteristic $p$, $R \rightarrow S$ be a surjection of graded $k$-algebras with $S$ artinian. Assume that either $R$ is a complete intersection or $R$ is Gorenstein and F-pure. If $d$ is the top socle degree of $S$, then the top socle degree  of $F^e_R (S)$ is at least $qd-(q-1)a(R)$.
\endproclaim

\demo{Proof} Write $S=R/J$, with $J\subset R$ an $\frak m$-primary ideal, where $\frak m$ is the unique homogeneous maximal ideal of $R$.

We first assume that $R$ is Gorenstein and F-pure. Let $\frak a$ be an $\frak m$-primary ideal of $R$, generated by a regular sequence, with $\frak a\subset J$. Let $g_1,\dots,g_s$, with $|g_1|\le \dots\le |g_s|$,  be elements in $R$ which represent a minimal generating set for $(\frak a\:\!\!  J)/\frak a$.  The hypothesis that $R$ is F-pure ensures that $g_1^q,\dots,g_s^q$ represents a minimal generating set for  $(g_1^q,\dots,g_s^q,\frak a^{[q]})/\frak a^{[q]}$; hence, the minimum degree among non-zero elements of  $(g_1^q,\dots,g_s^q,\frak a^{[q]})/\frak a^{[q]}$ is $q|g_1|$. Observation 1.4 yields that 
$$\tsize \operatorname{tsd} \frac RJ=\operatorname{socle\ degree}\frac R{\frak a }-|g_1|\quad\text{and}\quad \operatorname{tsd}\frac R{\frak a^{[q]}\:\!(g_1^q,\dots,g_s^q)}=\operatorname{socle\ degree}\frac R{\frak a ^{[q]}}-q|g_1|.$$
The $R$-module $R/\frak a$ has finite projective dimension so $$\tsize \operatorname{socle\ degree}\frac R{\frak a ^{[q]}}=
q\operatorname{socle\ degree}\frac R{\frak a }-(q-1)a(R).$$
Duality gives $J=\frak a\:\!(g_1,\dots,g_s)$. It follows that
$J^{[q]}\subseteq \frak a^{[q]}\:\!(g_1^q,\dots,g_s^q)$; and therefore, 
$$\align \tsize \operatorname{tsd}\frac R{J^{[q]}} \ge \operatorname{tsd}\frac R{\frak a^{[q]}\:\!(g_1^q,\dots,g_s^q)}&\tsize {}=\operatorname{socle\ degree}\frac R{\frak a ^{[q]}}-q|g_1|\\&\tsize {}=\operatorname{socle\ degree}\frac R{\frak a ^{[q]}}-q\left(\operatorname{socle\ degree}\frac R{\frak a }-\operatorname{tsd} \frac RJ\right)\\
&{}\tsize =q\operatorname{tsd} \frac RJ-(q-1)a(R) .\endalign $$

The proof is complete if $R$ is Gorenstein and F-pure. Throughout the rest of the argument, $R$ is a complete intersection.
We begin by reducing to the case where  $J$ is an irreducible ideal. Assume, for the time being, that the result has been established for irreducible ideals. 
Let $J=J_1\cap\dots \cap J_n$, with each $J_i$ irreducible.
Recall that   $\operatorname{tsd} R/J$  is the largest integer $d$ with $R_d\not\subseteq J$.
It follows that that the $\operatorname{tsd} R/J$ is equal to the maximum of the set $\{\operatorname{tsd} R/J_k\}$.   Fix a subscript $k$ with $\operatorname{tsd} R/J=\operatorname{tsd} R/J_k$. We know that $J^{[q]}\subseteq J_k^{[q]}$; and therefore,  
$$\tsize 
q\operatorname{tsd} \frac R J-(q-1)a(R)=q\operatorname{tsd} \frac R{J_k}-(q-1)a(R)\le \operatorname{tsd} \frac R{J_k^{[q]}}
\le
\operatorname{tsd} \frac R {J^{[q]}}. $$

Henceforth, the ideal $J$ is irreducible.  Write $R=P/C$, where $P$ is a polynomial ring and the ideal  $C$ is generated by the homogeneous regular sequence $f_1,\dots,f_c$. Let $I$ be the   pre-image of $J$ in $P$. In particular, $C\subseteq I$. The rings $R/J=P/I$ and $P/I^{[q]}$ are  Gorenstein, so Observation 1.4 gives 
$$\tsize \operatorname{socle\ degree} (\frac P{I^{[q]}})-M =\operatorname{tsd}(\frac P{I^{[q]}+C}),$$ where $M$ is the least degree among homogeneous non-zero elements of    $\frac {I^{[q]}\:\! C}{I^{[q]}}$. 
The $P$-module $P/I$ has finite projective dimension; so 
$$\tsize q\operatorname{socle\ degree} (\frac P{I})-(q-1)a(P) = \operatorname{socle\ degree} (\frac P{I^{[q]}}) .$$
 Recall the formula  $a(P/C)=a(P)+\sum_{i=1}^c|f_i|$.  The inequality 
$$\tsize q\operatorname{socle\ degree}(\frac P{I}) -(q-1) a(\frac PC) \le \operatorname{tsd} (\frac P{I^{[q]}+C})$$
is equivalent to the inequality 
$$\tsize M \le  (q-1)(|f_1|+\dots+|f_c|).$$
 We establish the most recent inequality. There exists an integer $t$, with $0\le t\le c(q-1)$, such that $C^{t}\not\subseteq I^{[q]}$; but $C^{t+1}\subseteq I^{[q]}$. Thus, some element   $f_1^{t_1}\cdots f_{c}^{t_c}$, with $\sum t_i=t$ and $0\le t_i\le q-1$ for all $i$,  of $C^{t}$ is a non-zero element of $\frac {I^{[q]}\:\! C}{I^{[q]}}$; and therefore,
$$\tsize M\le |f_1^{t_1}\cdots f_{c}^{t_c}|\le (q-1)(|f_1|+\dots+|f_c|).
\qed$$
\enddemo

The next result shows that we can get most of the way through the proof of Theorem 1.8 under the assumption that $R$ is Gorenstein and F-pure. Only the last step
(the result of Proposition 5.1) is still missing.
\proclaim{Proposition 6.2} Let $k$ be a field of positive characteristic $p$, $R \rightarrow S$ be a surjection of graded $k$-algebras with $R$ Gorenstein and $S$ artinian.  Assume, in addition, that $R$ is F-pure.
Assume that $d_1\le \dots\le d_{\ell}$ are the socle degrees of $S$ and the socle of $F_R^e(S)$ has the same dimension as the socle of $S$, with degrees of the generators given by $D_1\le D_2\le \cdots\le D_{\ell}$, with $$D_i=qd_i-(q-1)a(R),$$ for all $i$. 
Let $R=P/C$, and $S=R/IR$, with $P$ a polynomial ring, $I\subset P$. Then we have
$$
(C^{[q]} +I^{[q]})\:\! (C ^{[q]} \:\!C) =C + I^{[q]}.
$$
\endproclaim
\demo{Proof}
Let $A=C +(x_1, \ldots, x_d)$, where the images of $x_1, \ldots, x_d$ are a system of parameters in $R$. Let $K=A\:\!(I+C)$, so that we also have $I+C=A \:\!K$.
We have
$$(I ^{[q]} + C^{[q]} )\:\! (C^{[q]} \:\!C)= (A ^{[q]} \:\! K ^{[q]}) \:\!(C ^{[q]} \:\!C)=(A ^{[q]} \:\! (C ^{[q]} \:\!C)) \:\! K^{[q]}.
$$
We claim that $ A^{[q]} \:\!(C ^{[q]} \:\!C)= A ^{[q]} + C$. We see this by looking at the comparison map of resolutions induced by the projection $P/A^{[q]} \rightarrow P/(A^{[q]} +C)$. If $\Bbb F$ is the resolution of $P/C$, and ${\Bbb K}$ is the Koszul complex on $x_1, \ldots, x_d$, then the resolution of $P/A ^{[q]}$ is given by $\Bbb F^{[q]} \otimes {\Bbb K}^{[q]}$, the resolution of $P/(A ^{[q]} +C)$ is given by $\Bbb F \otimes {\Bbb K^{[q]}}$, and the comparison map between them is given by the comparison map $\Bbb F^{[q]} \rightarrow \Bbb F$, tensored with ${\Bbb K} ^{[q]}$. Thus, the last map is multiplication by an element of $P$ which represents the generator of $(C ^{[q]} \:\!C)/C ^{[q]}$.

It follows that $(I ^{[q]} + C ^{[q]} )\:\!( C ^{[q]} \:\!C) = (A ^{[q]} +C) \:\!K ^{[q]}$. It is clear that 
$$(A^{[q]}+C)\:\! K^{[q]}\supseteq I^{[q]}+C.$$ We next show that the rings defined by these ideals have the same socle degrees. 
Let $\delta$ and $\Delta$ be the socle degrees of the Gorenstein rings $\frac PA$ and $\frac P{A^{[q]}+C}$, respectively. The $P/C$-module $P/A$ has finite projective dimension; so $\Delta=q \delta-(q-1)a(\frac PC)$.
Let $g_1,\dots,g_s$ be elements of $K$ which represent a minimal generating set for $\frac KA$. 
Observation 1.4  gives that the socle degrees of $\frac P{I+C}$ are $\{\delta- |g_i|\}$.
So, our hypothesis tells us that the socle degrees of $\frac P{I^{[q]}+C}$ are
$$\tsize \{q(\delta- |g_i|)-(q-1)a(\frac PC)\}=\{\Delta- q|g_i|\}.$$
It is clear that $g_1^q,\dots,g_s^q$ represents a generating set for $\frac {K^{[q]}+C}{A^{[q]}+C}$. The hypothesis that $\frac PC$ is F-pure guarantees that $g_1^q,\dots,g_s^q$ is a minimal generating set of $\frac {K^{[q]}+C}{A^{[q]}+C}$. Observation 1.4  yields that the socle degrees of 
$\frac P{(A^{[q]}+C)\:\! K^{[q]}}$ are exactly the same as the socle degrees of $\frac P{I^{[q]}+C}$;  and therefore, the result of \cite{6} shows that   
$ I^{[q]} +C = (A ^{[q]} +C) \:\!K ^{[q]}= (I ^{[q]} + C ^{[q]} )\:\! (C ^{[q]} \:\!C)$.
\qed\enddemo

\enddocument